\documentclass{amsart}

\usepackage{ricky}

\begin{document}

\title{Quaternionic modular forms of any weight}

\author{Riccardo Brasca}

\email{\href{mailto:riccardo.brasca@gmail.com}{riccardo.brasca@gmail.com}}
\address{}

\curraddr{Max Planck Institute for Mathematics\\
Bonn\\
Germany}

\urladdr{\href{http://http://guests.mpim-bonn.mpg.de/ricky/}{http://http://guests.mpim-bonn.mpg.de/ricky/}}

\date{\today}



\begin{abstract}
In this work we construct an eigencurve for $p$-adic modular forms attached to an indefinite quaternion algebra over $\Q$. Our theory includes the definition, both as rules on test objects and sections of line bundle, of $p$-adic modular forms, convergent and overconvergent, of any $p$-adic weight. We prove that our modular forms can be put in analytic families over the weight space and we introduce the Hecke operators $\U$ and $\T_l$, that can also be put in families. We show that the $\U$-operator acts compactly on the space of overconvergent modular forms. We finally construct the eigencurve, a rigid analytic variety whose points correspond to systems of overconvergent eigenforms of finite slope with respect to the $\U$-operator.
\end{abstract}

\maketitle

\section*{Introduction}
Let $D$ be an indefinite quaternion algebra over $\Q$. Under certain hypotesis, the Jacquet-Langlands correspondence gives a map from the space of modular forms for $D$ to the space of modular forms for $\GL_2 / \Q$. In \cite{pay_thesis}, Kassaei developed the theory of $p$-adic modular forms (where $p$ is a prime that does not divide $\delta$, the discriminant of $D$), of integral weight, attached to $D$. The Jacquet-Langlands correspondence works also for $p$-adic modular forms, of level coprime with $p$. Since nowadays the theory of $p$-adic modular forms for $\GL_2/\Q$ is well developed, it is natural to seek for a similar theory also for modular forms relative to $D$. One of the most important result in the elliptic case is the existence of the eigencurve: a rigid space whose points correspond to systems of overconvergent eigenforms (of any $p$-adic weight) of finite slope with respect to the $\U_p$-operator. The final goal of this paper is the construction of the eigencurve for modular forms attached to $D$. We are going to use the general machinery developed by Buzzard in \cite{buzz_eigen}. We therefore need a weight space $\mc W$ and a notion of analytic family of modular forms over $\mc W$. Our weights are continuous characters $\Z_p^\ast \to K^\ast$, where $K$ is a finite extension of $\Q_p$, and we take as $\mc W$ the rigid analytic curve such that $\mc W(K) = \Hom_{\operatorname{cont}}(\Z_p^\ast,K^\ast)$. We define the space of modular forms, both convergent and overconvergent, with respect to $D$, level $N$ coprime with $p$, with coefficients in $K$ (where $K$ is big enough), and weight any $\chi \in \mc W(K)$. We then show that this modular forms can be put in families. We have the Hecke operators $\T_l$, for $l$ a prime not dividing $\delta N$, and the $\U$-operator, that acts completely continuously on the space of overconvergent modular forms. The main result of the paper (Theorem~\ref{teo: eigen}) is the following
\begin{teono}
There is a rigid space $\mc C \subseteq \mc W \times \m A^{1,\rig}_K$, called the \emph{eigencurve}, whose $L$-points, where $L$ is a finite extension of $K$, correspond naturally to systems of eigenvalues of overconvergent modular forms of finite slope with respect to the $\U$-operator, defined over $L$. Let $\pi_i$ be the projection to the $i$-th factor. If $x \in \mc C(L)$, let $\mc M_x$ be the set of overconvergent modular forms corresponding to $x$. Then the elements of $\mc M_x$ have weight $\pi_1(x) \in \mc W(L)$ and the $\U$-operator acts on $\mc M_x$ with eigenvalue $\pi_2(x)^{-1}$. We have that $\pi_1$ is, locally on $\mc W$ and on $\mc C$, finite and surjective.
\end{teono}

In \cite{shimura}, we have developed a theory of $p$-adic modular forms of any weight attached to a quaternion algebra over a totally real field $F$ different from $\Q$, following the approach of \cite{over}. With this paper we want to fill the gap left by the case $F=\Q$, using the same techniques. The main technical results of \cite{shimura} work also in the case $F = \Q$, but there are some differences, that we explain the following description of the paper.

In Section~\ref{sec: class mod forms}, we describe the basic objects of our work, following \cite{non-optimal}. We fix a prime $p > 3$ and an indefinite quaternion algebra $D$ over $\Q$ with discriminant $\delta$. We assume $p \nmid \delta$. Let $N$ be an integer with $N>4$ and coprime with $p$. We fix a suitable compact open subgroup $K(N) \subseteq D^\ast \otimes_{\Q} \m A^f$, where $\m A^f$ is the ring of finite adele of $\Q$. There is a smooth and proper PEL Shimura curve $\intN$ defined over $\Z_p$ that parametrizes abelian surfaces (more precisely `false elliptic curves') with $N$-level structure. We have also the curves $\intNpn{r}$ for any $r \geq 0$, analogue to the usual compactified modular curve of level $\Gamma_1(Np^r)$. Using the moduli problem solved by $\intNpn{r}$, we define the space of classical modular forms with respect to $D$, level $Np^r$, coefficients in $K$ (a finite extension of $\Q_p$), and weight $k \in \Z$. We empathize the approach of considering a modular form as a `rule defined on test objects' rather that a global section of a line bundle. This a big difference with \cite{shimura}, and we believe that it makes the theory more transparent (though we need to know the our modular forms actually are sections of some line bundle to construct the eigencurve). We take the rigid analytic viewpoint from the very beginning, so we will consider only modular forms over $K$.

In Section~\ref{sec: padic mod} we define $p$-adic modular forms of integral weights and we recall the theory of the canonical subgroup, following \cite{pay_thesis} and \cite{shimura}. In Section~\ref{sec: more can} we study the (Cartier dual of the) canonical subgroup in detail and we introduce the map $\dlog$, that will be crucial in the sequel. We obtain some of the results of \cite{over}, but our approach needs much less calculations and it is completely explicit. In particular, we prove that, over $\intNp$, the dual of the canonical subgroup admits a canonical point, and we study the image of this point under the map $\dlog$. We do this using a suitable Lubin-Tate group for $\Z_p$ instead of the formal multiplicative group. This idea comes from \cite{shimura}, where it was needed to take into account the action of a finite extension of $\Z_p$ on various objects. In our situation, this is not necessary, but it is in any case very useful, as we explain comparing our approach with the one of \cite{over}. In Section~\ref{sec: HT}, we explain the technical results about the map $\dlog$ that we will need.

In Section~\ref{sec: non int}, we recall the basic property of the weight space $\mc W$. Let $K$ be a finite extension of $\Q_p$ and let $\chi \colon \Z_p^\ast \to K^\ast$ be a continuous character. If $K$ contains all the $p^r$-th roots of unity, where $r$ depends on $\chi$, we define the space of $p$-adic modular forms with coefficients in $K$ and weight $\chi$, of any level. As mentioned above, our definition is quite simple and it is based on test objects. This makes the definition and the study of the diamond operators very clear. If $\chi(t)=t^k$, then, as usual, we obtain the space of modular forms of integral weight $k$. We then show that our modular forms are sections of a line bundle on a certain tubular neighborhood of the ordinary locus of the (rigid analytic) Shimura curve. We finally show that our modular forms can be put in analytic families over the weight space. This paves the way to the construction of the eigencurve via Buzzard's machinery.

In Section~\ref{sec: eigen} we finally construct the eigencurve. It remains to define the Hecke operators $\U$ and $\T_l$, where $l$ is a prime not dividing $\delta N$. Using test objects this is not difficult, but we need the approach via line bundles to prove that they are continuous and that $\U$ is completely continuous on the space of overconvergent modular forms. All the assumptions of \cite{buzz_eigen} are verified, so we obtain Theorem~\ref{teo: eigen}.
\subsection*{Acknowledgments}
I would like to thank Fabrizio Andreatta and Marco Seveso for their encouragements in writing this work. This paper was carried out when the author was a guest researcher at the Max Planck Institute for Mathematics. I would like to thank that institution for the hospitality and the excellent work conditions provided.

\section{Shimura curves and classical modular forms} \label{sec: class mod forms}
Let $p > 3$ be a rational prime, fixed from now on. Let $\val(\cdot)$ be the valuation of $\Q_p$ normalized in such a way that $\val(p)=1$. We will write $|\cdot|$ for the absolute value of $\Q_p$ such that $|x|= p^{-\val(x)}$ for all $x \in \Q_p$. We extend both $\val(\cdot)$ and $|\cdot|$ to the whole $\C_p$, that is the completion of a fixed algebraic closure of $\Q_p$. We will write $[\cdot] \colon \m F_p^\ast \to \Z_p$ for the Teichm\"uller character, and we set $[0] \colonequals 0$. All the schemes we are going to consider will be defined over $\Z_p$. When there is no possibility of confusion, we will denote base change with a subscript.

In this section we recall the theory of classical modular forms attached to a quaternion algebra, see \cite{non-optimal} and \cite{pay_thesis} for details. Let $D$ be a quaternion algebra over $\Q$ that is indefinite and with discriminant $\delta$ that satisfies $p \nmid \delta$. We assume that $D$ is a division algebra and we fix $\OD$, a maximal order of $D$. We fix identifications $D \otimes_{\Q} \R = \mat_2(\R)$ and $\OD \otimes_{\Z} \Z_l = \mat_2(\Z_l)$ for all primes $l \nmid \delta$. Let $x \mapsto x'$ be the canonical involution of $D$. We choose $t \in \OD$ such that $t^2 = -\delta$ and we define another involution on $D$ by $x \mapsto x^\ast \colonequals t^{-1}x't $. We write $\m A^f$ for the ring of \emph{finite} adele of $\Q$. Let $G / \Q$ be the reductive algebraic group such that its $R$-points, for any $\Q$-algebra $R$, are $G(R) = (D \otimes_{Q} R)^\ast$. We have that $G(\m A^f)$ is the restricted tensor product of the $D_l^\ast$'s, where $l$ is a prime and $D_l \colonequals D \otimes_{Q} \Q_l$. Let $\mc O_{D_l}$ be $\OD \otimes_{\Z} \Z_l$ and let $K$ be a compact open subgroup of $\prod_l \mc O_{D_l}^\ast$. We will assume that $K$ is `small enough', in particular the following conditions will suffice
\begin{itemize}
 \item $(N_K,\delta) = 1$, where $N_K$ (the \emph{level} of $K$) is the smallest positive integer $N$ such that $K$ contains all the elements of $\prod_l \mc O_{D_l}$ which are congruent to $1$ modulo $N$;
 \item $K \subseteq V_1(N)$ for some integer $N \geq 4$ such that $(N,\delta)=1$, where $V_1(N)$ is the subgroup of $\prod_l \mc O_{D_l}^\ast$ given by elements which are congruent to $\left(\begin{array}{cc} \ast & \ast \\ 0 & 1\end{array}\right)$ modulo $N$;
 \item $\det(K) = \widehat{\Z}^\ast$, where $\det \colon G(\m A^f) \to \m A^f$ is the reduced norm.
\end{itemize}
In particular we can take $K = V_1(N)$ with $N > 4$.

Let $T$ be any scheme over $\Z_p$. By a \emph{false elliptic curve} over $T$ we mean an abelian surface $\mc A \to T$ together with an embedding $i \colon \OD \hookrightarrow \End_T(\mc A)$. If $(N,\delta) = 1$, by a \emph{$N$-level structure} on a false elliptic curve $(\mc A,i)$ over $T$ we mean an isomorphism of left $\OD$-modules $\alpha \colon \mc A[N] \stackrel{\sim}{\to} (\OD \otimes_{\Z} \Z/N\Z)_T$. If $K$ is as above, a \emph{$K$-level structure} on $(\mc A,i)$ is a $N_K$-level structure defined up to right $K$-equivalence.

We now fix $K(N) \subseteq G(\m A^f)$ of level $N$ that satisfies the above assumptions and we assume that $(p, N) = 1$. It can be proved that, on any false elliptic curve $(\mc A,i)$ over $T$, there is a unique principal polarization whose corresponding Rosati involution on $\End(\mc A_x)$ (where $x$ is any geometric point of $T$) restricts to $^\ast$ on $\OD$. In particular, any false elliptic curve is canonically principally polarized. The functor $\Z_p\mbox{-schemes} \to \mathbf{set}$ that sends $T$ to the set of isomorphism classes of false elliptic curves with $K$-level structure is representable by a geometrically connected scheme $\intN$ defined over $\Z_p$. We have that $\intN \to \Sp(\Z_p)$ is smooth, proper, and of relative dimension $1$. The universal object will be denoted $\intAN \to \intN$.

Let $\mc C$ be a pseudo-abelian category and let $X$ be an object of $\mc C$ with an action of $\OD \otimes_{\Z} \Z_p$. Choosing a non-trivial idempotent $e \in \OD \otimes_{\Z} \Z_p \cong \mat_2(\OP)$ such that $e^\ast = e$, we obtain a decomposition
\[
X = X^1 \oplus X^2,
\]
where each $X^i$ has an action of $\Z_p$. Note that $X^1$ and $X^2$ are isomorphic. This notation will be used throughout the paper.

Let $n \geq 0$ be an integer and let $K(N,p^n)$ be the subset of $K(N)$ given by the elements whose component at $p$ is upper triangular modulo $p^n$. There is a regular scheme $\intNsubpn{n}$ that represents the functor $\Z_p\mbox{-schemes} \to \mathbf{set}$ that sends $T$ to the set of isomorphism classes of quadruples $(\mc A,i,\alpha, \mc C)$ where $(\mc A,i)$ is a false elliptic curve over $T$, $\alpha$ is a $K$-level structure on $(\mc A,i)$ and $\mc C$ is a finite and flat subgroup of order $p^n$ of $\mc A[p^n]^1$. We have that $\intNsubpn{n} \to \Sp(\Z_p)$ is proper and of relative dimension $1$, generically smooth but not smooth. Indeed, the special fiber is a normal crossing divisor with two components, each one isomorphic to $\intN_{\m F_p}$. The two components intersect at the supersingular points (see below). The universal object will be denoted $\intANsubpn{n} \to \intNsubpn{n}$.
\begin{rmk} \label{rmk: can sub}
Any $\mc C$ as in the above moduli problem gives $\mc D$, a finite and flat subgroup of order $p^n$ of $\mc A[p^n]^2$, and we have that $\mc C \oplus \mc D$ is a finite and flat subgroup of $\mc A[p^n]$ of order $p^{2n}$, stable under the action of $\OD$. All such subgroups of $\mc A[p^n]$ arise in this way. The canonical principal polarization of $\mc A$ induces isomorphism $\mc A[p^n]^1 \cong (\mc A[p^n]^2)^\vee$ (Cartier dual). This isomorphism takes $\mc C \hookrightarrow \mc A[p^n]^1$ to $(\mc A[p^n]^2 / \mc D)^\vee \hookrightarrow \mc D^\vee$.
\end{rmk}
Let $K(Np^n)$ be the subset of $K(N,p^n)$ given by the elements whose component at $p$ is congruent to $\left(\begin{array}{cc} \ast & \ast \\ 0 & 1\end{array}\right)$ modulo $p^n$. There is a regular scheme $\intNpn{n}$ that represents the functor $\Z_p\mbox{-schemes} \to \mathbf{set}$ that sends $T$ to the set of isomorphism classes of quadruples $(\mc A,i,\alpha, P)$ where $(\mc A,i,\alpha)$ is as above and $P$ is a point of exact order $p^n$ of $\mc A[p^n]^1$. We have that $\intNpn{n} \to \Sp(\Z_p)$ is proper, generically smooth, and of relative dimension $1$. The universal object will be denoted $\intANpn{n} \to \intNpn{n}$.

The curves $\intN$, $\intNsubpn{n}$, and $\intNpn{n}$ are the Shimura curves we are interested in. There are several morphisms between them, given the by obvious natural transformations of functors. Let $K$ be any one of $K(N)$, $K(N,p^n)$ or $K(Np^n)$ and let $\mc M$ be the corresponding curve. If $e \colon \mc M \to \mc A$ is the zero section of the universal object, we define
\[
\underline \omega \colonequals \underline \omega_K \colonequals \left(e^\ast \Omega^1_{\mc A / \mc M} \right)^1,
\]
that is a locally free sheaf of rank $1$ on $\mc M$. Given $\Sp(R) \to \mc M$, where $R$ is a $\Z_p$-algebra, the pullback of $\underline \omega$ to $\Sp(R)$ will be denoted $\underline \omega_R$ or $\underline \omega_{\mc A/R}$, where $\mc A$ is the pullback of the universal object to $\Sp(R)$.

Let $V$ be a finite extension of $\Z_p$, with field of fraction $K$ (we will make several assumptions about $V$ during the paper). From now on we will usually work over $V$ or over $K$. To simplify the notation, we will omit the subscripts $_V$ and $_K$. We will also be interested in formal schemes and we will use the following convention: algebraic objects defined over $V$ will be denoted using Italics letter, like $\intN$. The completion along the subscheme defined by $p=0$ will be denoted using the corresponding Gothic letter, like $\formN$.
\begin{defi} \label{defi:test object}
Let $r \geq 0$ be an integer. A (classical) test object of level $K(Np^r)$ is a couple $((\mc A/R, i, \alpha, P), \omega)$, where:
\begin{itemize}
 \item $R$ is a $V$-algebra;
 \item $(\mc A / R, i, \alpha, P)$ is an object of the moduli problem of level $K(Np^r)$ (we will usually omit $P$ from the notation if $r=0$), with $\mc A$ defined over $R$ and such that $\underline \omega_R$ is free;
 \item $\omega$ is a basis of $\underline \omega_R$.
\end{itemize}
\end{defi}
\begin{defi} \label{defi: mod form}
Let $k$ be an integer. A (classical) \emph{modular form with respect to $D$, level $K(Np^r)$ and weight $k$, with coefficients in $K$}, is a rule that assigns to every test object $T = ((\mc A/R, i, \alpha, P), \omega)$, where $T$ is of level $K(Np^r)$, an element $f(T) \in R_K$ such that
\begin{itemize}
 \item $f(T)$ depends only on the isomorphism class of $T$;
 \item if $\varphi \colon R \to R'$ is a morphism of $V$-algebras and we denote with $T'$ the base change of $T$ to $R'$, we have $f(T') = \varphi(f(T))$;
 \item for any $\lambda \in R^\ast$, we have
\[
f(((\mc A/R, i, \alpha, P), \lambda\omega)) = \lambda^{-k} f(T).
\]
\end{itemize}
We write $S^D(K,K(Np^r),k)$ for the space of such modular forms.
\end{defi}
\begin{prop} \label{prop: mod form sheaf}
There is an invertible sheaf $\Omega^k(Np^r)$ on $\formNpn{r}^{\rig}$ such that
\[
S^D(K,K(Np^r),k) = \Homol^0(\formNpn{r}^{\rig}, \Omega^k(Np^r)).
\]
\end{prop}
\begin{proof}
We take for $\Omega^k(Np^r)$ the rigidification of formal completion of $\underline \omega^{\otimes k}$.
\end{proof}
\begin{rmk} \label{rmk: trivial}
The above proposition is trivial, and we can state a much stronger result (namely that we have an invertible sheaf over $\intNpn{r}$ that is a model of $\Omega^k(Np^r)$), but later on we will consider $p$-adic modular forms of any weight, and the above result has a direct analogue, while the stronger version has not.
\end{rmk}
Definition \ref{defi: mod form} is given only for modular forms over $K$, since this is the case we will be mainly interested in, but using Proposition \ref{prop: mod form sheaf}, we have the notion of classical modular form defined over any $\Z_p$-algebra, as usual.

\section{\texorpdfstring{$p$-adic modular forms of integral weight}{p-adic modular forms of integral weight}} \label{sec: padic mod}
We now analyze more closely what happens at level $K(N)$. There is a Kodaira-Spencer isomorphism $\underline \omega^{\otimes 2} \cong \Omega^1_{\intN/\Z_p}$ and we have $\deg(\underline \omega) = g-1$, where $g$ is the genus of any geometric fiber of $\intN$. Let $(\mc A,i,\alpha)$ be an object of the moduli problem, with $\mc A$ defined over a $V$-algebra $R$. We have that $\mc A[p^\infty]^1$ is a $p$-divisible group of dimension $1$ and height $2$. Let $\mathfrak A$ be the $p$-adic completion of $\mc A$ and let $\widehat{\mc A}$ be the completion of $\mathfrak A$ along the zero section. Then $\widehat{\mc A}^1$ is a formal group of dimension $1$. If $\widehat{\mc A}^1$ is a Barsotti-Tate group, its height is either $1$ or $2$. We say that $(\mc A,i, \alpha)$, or simply $\mc A$, is \emph{ordinary} if $\widehat{\mc A}^1$ has height $1$. Otherwise we say that $\mc A$ is \emph{supersingular}. We have that $\mc A$ is ordinary (supersingular) if and only if $\mc A_{\overline{\m F}_p}$ is ordinary (supersingular) in the usual sense. If $\underline \omega_R$ is a free $R$-module, we have that $\widehat{\mc A}^1$ is a formal group law, and there is a coordinate $x$ on $\widehat{\mc A}^1$ such that the multiplication by $p$ has the form
\[
[p](x)=p x+ a x^p + \sum_{j=2}^{\infty}c_jx^{j(p-1)+1},
\]
where $a$, $c_j$ are in $R$ and $c_j \in pR$ unless $j \equiv 1 \bmod{p}$. If we assume that $p=0$ in $R$, the various $a$'s glue together to define $H$, a modular form of level $K(N)$ and weight $p-1$, defined over $\m F_p$, that is called the \emph{Hasse invariant}. If $W=\Sp(R)$ is an open affine of $\intN$ and we denote with $\omega$ the differential dual to the coordinate $x$ defined above, we have $H_{|W}=a \omega^{\otimes p-1}$. We have that $H$ vanishes precisely at the supersingular points of $\intN_{\overline{\m F}_p}$ and has simple zeros here (so the number of supersingular points is $(p-1)(g-1)$). Moreover, supersingular points actually exist (hence $g \geq 2$). The modular form $H$ can be lifted, non-canonically, to a modular form $E_{p-1}$ defined over $\Z_p$. We fix once and for all such a lift, everything will not depend on this choice. Over $\Sp(R)$ we can write ${E_{p-1}}_{|\Sp(R)} = E \omega^{p-1}$ for some $E \in R$. We have $a \equiv E \bmod{p}$.

Let $0 \leq w < 1$ be a rational number, and let us assume that $V$ contains an element of valuation $w$, denoted $p^w$. We define
\[
\intN(w)\colonequals \Sp_{\intN}(\Sym(\underline \omega^{\otimes p-1})/(E_{p-1}-p^w)).
\]
Note that $\intN(w)$ has a natural moduli interpretation. It classifies couples $((\mc A, i, \alpha), Y)$, where $(\mc A/R,i,\alpha)$ is false elliptic curve of level $K(N)$ and $Y$ is a global section of $\underline \omega_{\mc A/R}^{\otimes 1-p}$ such that $Y E_{p-1}= p^w$.

Taking the rigidification of the map $\formN(w)\to \formN$ we get the immersion $\formN^{\rig}(w) \hookrightarrow \formN^{\rig}$, where $\formN^{\rig}(w)$ is the affinoid subdomain of $\formN^{\rig}$ defined by Coleman in \cite{cole_class}, relative to $E_{p-1}$. We call $\formN(0)^{\rig}$ the \emph{ordinary locus}, it is an affinoid subdomain of $\formN^{\rig}$: its complement is a finite union of open discs of radius $1$, called the \emph{supersingular discs}.

Let $r \geq 1$ be an integer and assume that $0 \leq w < \frac{1}{p^{r-2}(p+1)}$. Let $\mc A$ be an object of the moduli problem of $\intN(w)$, defined over $R$. In this case $\mc A[p^r]^1$ has a canonical finite and flat subgroup $\mc C_r$ of order $p^r$. This gives a morphism $\formN(w) \to \formNsubpn{r}$. Its rigidification is a section, defined over $\formN(w)^{\rig}$, of the morphism $\formNsubpn{r}^{\rig} \to \formN^{\rig}$. We define $\formNpn{r}(w)^{\rig}$ as the inverse image of $\formN(w)^{\rig}$ with respect to the map $\formNsubpn{r}^{\rig} \to \formN^{\rig}$. It is an affinoid subdomain of $\formNpn{r}^{\rig}$ with a finite and \'etale map to $\formN(w)^{\rig}$. Taking the normalization of $\formN(w)$ in $\formNpn{r}(w)^{\rig}$, that is a finite extension of its generic fiber, we obtain $\formNpn{r}(w)$, a formal model of $\formNpn{r}(w)^{\rig}$.

If we work `at level $K(Np^r)$', we will tacitly assume that $w < \frac{1}{p^{r-2}(p+1)}$, so we have the canonical subgroup of any object of the moduli problem of $\intNpn{r}(w)$. In particular, all our results will holds if $w$ is small enough.
\begin{notation}
Let $R$ be a $p$-adically complete and flat $V$-algebra. We will say that $R$ is \emph{normal} if $R$ is integrally closed in $R_K$. We extend this notion to formal schemes in the obvious way.
\end{notation}
\begin{prop}
We have that $\formNpn{r}(w)$ is a flat and normal formal scheme over $\Spf(V)$.
\end{prop}
\begin{proof}
In the case $r=0$ this is \cite[Proposition~9.5]{pay_totally}. If $r \geq 1$ it is a trivial consequence of the definition.
\end{proof}
\begin{prop} \label{prop: mod Npr}
Let $R$ be a normal, flat, and $p$-adically complete $V$-algebra. There is a natural bijection between $\formNpn{r}(w)(R)$ and the set of isomorphism classes of couples $((\mc A/R, i, \alpha, P), Y)$, where:
\begin{itemize}
 \item $(\mc A/R, i,\alpha, P)$ is an object of the moduli problem, with $\mc A$ defined over $R$, of $\intNpn{r}$ and $P$ generates the canonical subgroup of $\mc A[p^r]^1$;
 \item $Y$ is a section of $\underline \omega_{\mc A/R}^{\otimes 1-p}$ that satisfies $YE_{p-1}=p^w$.
\end{itemize}
\end{prop}
\begin{proof}
This is \cite[Proposition~3.2]{shimura}.
\end{proof}
The above proposition motivates the following
\begin{defi} \label{defi:test object p-adic}
Let $r \geq 0$ be an integer. A ($p$-adic) test object of level $K(Np^r)$ and growth condition $w$ is a couple $((\mc A/R, i, \alpha, P, \omega), Y)$, where:
\begin{itemize}
 \item $R$ is a $p$-adically complete, flat, and \emph{normal} $V$-algebra;
 \item $((\mc A / R, i, \alpha, P), \omega)$ is a classical test object of level $K(Np^r)$;
 \item $Y$ is a global section of $\underline \omega_{\mc A/R}^{\otimes 1-p}$ such that $YE_{p-1} = p^w$.
\end{itemize}
\end{defi}
\begin{defi} \label{defi: p-adic mod form}
Let $k$ be an integer. A \emph{$p$-adic modular form with respect to $D$, level $K(Np^r)$ and weight $k$, with coefficients in $K$ and growth condition $w$}, is a rule that assigns to every test object $T = ((\mc A/R, i, \alpha, P, \omega), Y)$, where $T$ is of level $K(Np^r)$ and growth condition $w$, an element $f(T) \in R_K$ such that
\begin{itemize}
 \item $f(T)$ depends only on the isomorphism class of $T$;
 \item if $\varphi \colon R \to R'$ is a morphism of $p$-adically complete, flat, and normal $V$-algebras and we denote with $T'$ the base change of $T$ to $R'$, we have $f(T') = \varphi(f(T))$;
 \item for any $\lambda \in R^\ast$, we have
\[
f((\mc A/R, i, \alpha, P, \lambda\omega), Y) = \lambda^{-k} f(T).
\]
\end{itemize}
We write $S^D(K,w,K(Np^r),k)$ for the space of such modular forms.
\end{defi}
If $h \leq r$ are integers, the natural injective morphism $S^D(K,w,K(Np^h),k) \to S^D(K,w,K(Np^r),k)$ allows us to consider $S^D(K,w,K(Np^h),k)$ as a subset of $S^D(K,w,K(Np^r),k)$.
\begin{prop} \label{prop: p-adic mod form sheaf}
There is an invertible sheaf $\Omega^k(Np^r)(w)$ on $\formNpn{r}(w)^{\rig}$ such that
\[
S^D(K,w,K(Np^r),k) = \Homol^0(\formNpn{r}(w)^{\rig}, \Omega^k(Np^r)(w)).
\]
\end{prop}
\begin{proof}
We take for $\Omega^k(Np^r)(w)$ the restriction to $\formNpn{r}(w)^{\rig}$ of $\Omega^k(Np^r)$.
\end{proof}
This proposition shows that the natural morphism $S^D(K,w',K(Np^r),k) \to S^D(K,w,K(Np^r),k)$, where $w' \geq w$ is a rational number that satisfies the same assumptions as $w$ does, is injective.
\begin{rmk} \label{rmk: trivial 2}
As in Remark~\ref{rmk: trivial}, we can state a stronger result than in the above proposition, namely, that we have a formal model of $\Omega^k(Np^r)(w)$. But, since we do not have an integral model of $\formNpn{r}(w)$ that satisfies a reasonable moduli property, we do not have an integral model of $\Omega^k(Np^r)(w)$. Later on, we will not even have the formal model of the sheaves we are interested in.
\end{rmk}
By definition, elements of $S^D(K, K(Np^r), k)$ correspond to sections of $\Omega^k(Np^r)$ over $\formNpn{r}^{\rig}$, while elements of $S^D(K,w,K(Np^r),k)$ correspond to sections over $\formNpn{r}(w)^{\rig}$. Elements of $S^D(K,0,K(Np^r),k)$ are called \emph{convergent} modular forms, while the elements of $S^D(K,w,K(Np^r),k)$, for $w>0$, are called \emph{overconvergent} modular forms.

\section{\texorpdfstring{More on the canonical subgroup and the map $\dlog$}{More on the canonical subgroup and the map dlog}} \label{sec: more can}
There is a primitive $p$-th root of unity $\zeta_p$ in $V$ if and only $V$ contains $\xi$, a non trivial root of $x^p - px$, that is, a $(p-1)$-th root of $-p$. This is an exercise, but it also follows from Oort-Tate theory: by \cite[pages 8-10]{oort}, there is an isomorphism of group schemes $\Sp(V[x]/(x^p-1)) \cong \Sp(V[x]/(x^p-pw_{p-1}x))$, where $w_{p-1} \in \Z_p^\ast$ is a universal constant. In particular, $\zeta_p$ gives a $(p-1)$-th root of $w_{p-1} p$. We have that the reduction modulo $p$ of $w_{p-1}$ is congruent to $(p-1)!$, so $-w_{p-1}$ admits a $(p-1)$-th root in $\Z_p$ and we get the desired $(p-1)$-th root of $-p$. The point of this discussion is that, to simplify the calculations, we prefer to take Oort-Tate theory as \emph{starting point} rather than using it to obtain technical results. In this way our calculations will be completely trivial. From this point of view, $\mu_{p,V}$ \emph{is} $\Sp(V[x]/(x^p+px))$. From now on, when we work at level $K(Np)$, we will tacitly assume that $V$ contains a fixed $(p-1)$-th root of $-p$, denoted $(-p)^{1/(p-1)}$ (in particular, $V$ contains also $\zeta_p$). Later on, we will explain a similar assumption for higher levels, equivalent to the fact that $V$ contain a primitive $p^r$-root of unity.

At level $K(Np)$, the canonical subgroup is given, locally on $\intN(w)$, by $\mc C = \mc C_1 = \Sp(R[x]/(x^p+\frac{p}{E}x))$, where ${E_{p-1}}_{|\Sp(R)} = E \omega^{p-1}$ as in Section~\ref{sec: padic mod} (the fact that the canonical subgroup arise naturally in this form is the main reason why it is so convenient to use Oort-Tate theory from the very beginning).
\begin{prop} \label{prop: E exist}
There is $E_1 \in S^D(K,w,K(Np),1)$ such that $E_1^{q-1}=E_{q-1}$.
\end{prop}
\begin{proof}
Let $\Sp(S)$ be the base change of $\Sp(R)$ to $\intNp(w)$. By the moduli property of $\intNp(w)$, we have a canonical $S$-point of $\mc C$, that gives (using $(-p)^{1/(p-1)}$) a $(p-1)$-th root of $E$ in $S$. These roots glue together to define the required $E_1$.
\end{proof}
Let $r \geq 1$ be an integer. Fix an open affine $\Spf(R)$ of $\formN(w)$ such that $\underline \omega$ is free, generated by $\omega$ as above. Let $\Spf(S_r)$ be the pullback of $\Spf(R)$ to $\formNpn{r}(w)$ under $\formNpn{r}(w) \to \formN(w)$. Let $\mc A$ be any object of the moduli problem of $\formN(w)$. We know that $\mc A[p^r]^1$ admits a canonical subgroup $\mc C_r$ and that the canonical point of $\mc A_{S_r}[p^r]^1$ generates $\mc C_r$. We want to construct a canonical $S_r$-point of $(\mc C_r)^\vee$ that is a generator of $(\mc C_r)^\vee(S_r)$ as $\Z/p^r\Z$-module. All our schemes are normal, so it is enough to define a morphism $\mc C_{r,S_{r,K}} \to \mu_{p^r,S_{r,K}}$. But $\mc C_{r,S_{r,K}}$ is a constant group scheme, with $\Z/p^r\Z$ as associated abstract group, and we have $P$, a canonical generator of $\mc C_{r,S_{r,K}}$. So it is enough to define the image of $P$, or, that is the same, to choose a point of $\mu_{p^r,S_{r,K}}$. In particular it is natural to assume that $V$ contains $\zeta_{p^r}$, a fixed primitive $p^r$-th root of unity, and to take it as the image of $P$. But recall that we want to use Oort-Tate theory as building block, so for us $\mu_{p,V}$ is $\Sp(V[x]/(x^p+px))$ rather than $\Sp(V[x]/(x^p-1))$. For this reason, we take the following slightly different approach.

Let $\mc{LT}$ be the \emph{basic Lubin-Tate group associated to $p$}, as in \cite[Chapter 8]{lang_cycl}, over $V$. It is a formal group law on $V[[x]]$ that has an action of $\Z_p$ and the multiplication by $p$ is given by $[p](x)=x^p+px$. If $z \in \m F_p$, the action of $[z]$ sends $x$ to $[z]x$. All Lubin-Tate groups associated to $p$ are strictly isomorphic, so $\mc{LT} \cong \widehat{\m G}_{\operatorname m, V}$. In particular, if $G$ is a finite and flat group scheme over $S_r$, giving a morphism $G \to \mc{LT}_{S_r}$ is the same as giving a morphism $G \to \widehat{\m G}_{\operatorname m, S_r}$. Repeating the above discussion, we see that to define a point of $(\mc C_r)^\vee(S_r)$ it is enough to fix a point of $\mc{LT}$ of order exactly $p^r$. For this reason we fix a sequence $\set{\xi_r}_{r \geq 1}$ of $\C_p$-points of $\mc{LT}$ such that $\xi_r$ has order exactly $p^r$ and $p\xi_{r+1} = \xi_r$ (recall that $\mc{LT}$ is denoted additively). By definition, the points of order $p$ of $\mc{LT}$ are the points of $\Sp(V[x]/(x^p+px))$, so it is natural to require that $\xi_1$ is our fixed $(-p)^{1/(p-1)}$. When we will work at level $K(Np^r)$, we will tacitly assume that $\xi_r \in V$, so we have $\gamma_r$, a canonical $S_r$-point of $(\mc C_r)^\vee$ of order $p^r$. For various $r$'s, these morphisms are compatible, and if we let $r \to \infty$ (and we take $w=0$) we have a morphism, defined over $\C_p$
\[
\mc C_\infty \to \mc{LT}[p^\infty],
\]
where $\mc C_\infty$ is the $p$-divisible group defined by the $\mc C_r$'s.

We can study more closely the case $r=1$. Recall that the canonical subgroup is $\mc C = \Sp(R[x]/(x^p + \frac{p}{E}))$, where $E \in R$ satisfies ${E_{p-1}}_{|\Spf(R)} = E \omega^{\otimes p-1}$. Let $S \colonequals S_1$. Giving $\gamma = \gamma_1$, our canonical $S$-point of $\mc C_S^\vee$, is the same as giving the morphism $S[x](x^p+px) \to S[x]/(x^p + \frac{p}{E})$, that amounts to choose a $(p-1)$-th root of $E$ in $S$. But by definition we have a $(p-1)$-th root of $-\frac{p}{E}$ in $S$, and using our fixed $\xi_1$ we get $E^{1/(p-1)} \in S$. At the end we see that our morphism $\gamma \in \mc C^\vee(S)$ corresponds to
\begin{gather*}
\mc C_S = \Sp(S[x]/(x^p+\frac{p}{E})) \to \mc {LT}_S = \Spf(S[[x]]) \\
E^{1/(p-1)} x\mapsfrom x
\end{gather*}
In this way the explicit description of $\gamma$ is very simple, and the relation between $\gamma$ and $E_1$ is completely transparent.

It follows by the explicit description of the canonical subgroup that the module of invariant differential of $\mc C$ is generated by $\di(x)$, with the unique relation $\frac{p}{E}\di(x)=0$, so we have $\underline \omega_{\mc C} \cong R/\frac{p}{E}R\di(x)$. We can be even more precise: let $h \colon \mc C \to \widehat{\mc A}[p]^1$ be the natural map. We have fixed above a coordinate $x$ of $\widehat{\mc A}^1$, and we have denoted with $\omega$ the differential dual to $x$. Then, in $\Omega^1_{\mc C/R}$, we have $h^\ast(\omega) = \frac{\di(x)}{1-E}$. In particular, under the isomorphism $\underline \omega_{\mc C} \cong R/\frac{p}{E}R\di(x)$, we have $h^\ast(\omega) = \di(x)$. Using Oort-Tate theory, we can give an explicit formula for the comultiplication of $\mc C$. It is given by
\[
c(x)= x \otimes 1 + 1\otimes x + \frac{1}{p-1}\sum_{i=1}^{p-1}E\frac{w_{p-1}}{w_i w_{p-1}}x^i \otimes x^{p-i},
\]
where $w_i \in \Z_p$ are the universal constants of \cite{oort}. Moreover, the counit of $\mc C$ is given by $x \mapsto 0$. Having fixed $\zeta_p$, a primitive $p$-th root of unity in $V$, one can also write down a morphism $\mc C_S \to \Sp(S[x]/(x^p-1))$ explicitly (clearly, to obtain the morphism corresponding to $\gamma$, we need to take as $\zeta_p$ the root of unity corresponding to $\xi_1$ under the isomorphism $R[x]/(x^p-1) \cong R[x]/(x^p+p)$). This is the approach taken in \cite{over}, and we now show that our point of view makes the calculations much easier. Let $\alpha \in S$ be $(p-1)$-th root of $-\frac{p}{E}$ given by the moduli problem (so we have $\xi_1 = \alpha E^{1/(p-1)}$). Then the map
\begin{gather*}
S[x]/(x^p-1) \to S[x]/(x^p+\frac{p}{E}x) \\
x \mapsto 1 + \frac{1}{1-p}\sum_{i=1}^{p-1} \frac{(-1)^i x^i}{w_i \alpha^{i}}\left(\sum_{z \in \m F_p^\ast}[x]^{-1}\zeta_p^z \right)^i
\end{gather*}
gives the required morphism $\mc C_S \to \Sp(S[x]/(x^p-1))$. We see that here the formula is much more complicated than the one for $\gamma$, in particular, Proposition~\ref{prop: congr E} below, that is easy with our approach, requires non trivial calculations in \cite{over}.

We continue to use the above notations. Let $\eta = \Sp(\m K)$ be a generic geometric point of $\Sp(R)$, we write $\mc G$ for $\pi_1(\Sp(R_K), \eta_K)$. We denote with $\overline R$ the direct limit of all $R$-algebras $T \subseteq \m K$ which are normal and such that $T_K$ is finite and \'etale over $R_K$. Let $G$ be a finite and flat group scheme over $R$, killed by $p^n$. Taking the pullback of the invariant differential $\di(x)$ of $\mc{LT}$, we get the map
\[
\dlog_{G,T} \colon G^\vee(T) \to \underline \omega_G \otimes_R T/p^nT,
\]
where $T$ is as above. We are interested in the map $\dlog_G \colon G^\vee(\overline R) \to \underline \omega_G \otimes_R \overline R/p^n \overline R$ obtained by direct limit. Taking the inverse limit over $n$, we have the map
\[
\dlog_{\mc A} \colon \mc \T_p((A[p^\infty]^1)^\vee) \otimes_{\Z_p} \widehat{\overline R} \to \underline \omega_{\mc A/R} \otimes_R \widehat{\overline R},
\]
where $\T_p((A[p^\infty]^1)^\vee) = ((A[p^\infty]^1)^\vee)(\overline R)$ is the Tate module of $\mc A[p^\infty]^1$.
\begin{prop} \label{prop: congr E}
We have
\[
\dlog_{\mc A[p]^1,S}(\gamma) \equiv E_1 \bmod{p^{1-w}}.
\]
\end{prop}
\begin{proof}
Clearly we have $\gamma^\ast(\di(x))= E^{1/(p-1)}\di(x)$. We can write $\underline \omega_{\mc C} = \omega R/\frac{p}{E}R$ and $\underline \omega_{\mc A} = \omega R$. But by definition of $\intN(w)$, we have $\frac{p}{E} \in p^{1-w} R$, so the natural map $\underline \omega_{\mc A} \to \underline \omega_{\mc C}$ is an isomorphism modulo $p^{1-w}$. By definition we have ${E_1}_{|\Spf(R)} = E^{1/(p-1)} \omega \in \Homol^0(\formNp(w), \underline \omega)$ so we see that $\gamma^\ast(\di(x))\equiv {E_1}_{|\Spf(S)} \bmod{p^{1-w}}$.
\end{proof}
From now on we will omit $(\overline R)$ in the notation, it should be clear from the context whether we are talking about the group scheme or about the group of points. We also write $\overline R_z$ for $\overline R/p^z \overline R$ (and similarly for other objects).

\section{The Hodge-Tate sequence} \label{sec: HT}
We recall in this section the main technical results of \cite{shimura} (see also \cite{over}). The proof are exactly the same, so we will omit them. Let $\Spf(R)$ be an open affine of $\formN(w)$, and let $\mc A$ be the corresponding false elliptic curve defined over $R$. We assume that $w < \frac{p}{p+1}$.
\begin{defi} \label{defi: HT seq}
The \emph{Hodge-Tate sequence of $\mc A$} is the following sequence of $\widehat{\overline R}$-modules with semilinear action of $\mc G$:
\[
0 \to \underline \omega_{\mc A^\vee/R}^\ast \otimes_R \widehat{\overline R}(1) \stackrel{a_{\mc A}}{\longrightarrow} \T_p((\mc A[\varpi^\infty]^1)^\vee) \otimes_{\Z_p} \widehat{\overline R} \stackrel{\dlog_{\mc A}}{\longrightarrow} \underline \omega_{\mc A/R} \otimes_R \widehat{\overline R} \to 0.
\]
where $^\ast$ means `dual module', $a_{\mc A}$ is defined using the $\dlog$-maps of $(\mc A[p^\infty]^1)^\vee$, and $(1)$ is the Tate-twist.
\end{defi}
\begin{teo} \label{teo: HT}
The Hodge-Tate sequence is a complex and the map $a_{\mc A}$ is injective. Moreover, the homology of the sequence is killed by $p^v$, where $v \colonequals \frac{w}{p-1}$. The Hodge-Tate sequence is exact if and only if $\mc A$ is ordinary. Furthermore, $\ker(\dlog_{\mc A})$ and $\im(\dlog_{\mc A})$ are free $\widehat{\overline R}$-modules of rank $1$. The latter is generated by $\dlog_{\mc A}((\gamma)_r)$. Finally, taking $\Gal(\overline R_K/S_K)$-invariants, everything remains true also over $S$ (i.e. at level $K(Np)$), but we cannot descend to $R$ (i.e. at level $K(N)$).
\end{teo}
\begin{coro} \label{coro: glob free}
We have that $\underline \omega_{K(Np^r)}^{\rig}$ is \emph{globally} free if and only if $r \geq 1$. In this case, we have $\underline \omega_{K(Np^r)}^{\rig} = \mc O_{\formNpn{r}(w)^{\rig}} E_1$.
\end{coro}
This corollary is the reason why the theory is quite different if $r \geq 1$ or not.

Let $\mc D_r \subseteq (\mc A[p^r]^1)^\vee $ be the orthogonal of $\mc C_r$. Let $w < \frac{1}{p^{r-2}(p+1)}$ as usual (we need to assume $w \leq 1/p$ if $r=1$). Then the kernel of the map $\dlog_{\mc A[p^r]^1, \overline R} \colon (\mc A[p^r]^1)^\vee \to \underline \omega_{\mc A/R} \otimes_R \overline R_r$ is $\mc D_r$.
\begin{prop} \label{prop: HT}
We have a commutative diagram of $\mc G$-modules, with exact rows and vertical isomorphisms,
\[
\xymatrix{
0 \ar[r] & \ker(\dlog_{\mc A})_{r-v} \ar[r] \ar[d]^\wr & \T_p((\mc A[p^\infty]^1)^\vee) \otimes \overline R_{r-v} \ar[r] \ar@{=}[d] & \im(\dlog_{\mc A})_{r-v} \ar[r] \ar[d]^\wr & 0\\
0 \ar[r] & \mc D_r \otimes \overline R_{r-v} \ar[r] & (\mc A[p^r]^1)^\vee \otimes \overline R_{r-v} \ar[r] & \mc C_r^\vee \otimes \overline R_{r-v} \ar[r] & 0
}
\]
\end{prop}

\section{Modular forms of non-integral weight} \label{sec: non int}
There is a rigid analytic (group) variety $\mc W$ over $\Q_p$ such that, for any affinoid algebra $A$, we have $\mc W(A) = \Hom_{\operatorname{cont}}(\Z_p^\ast,A^\ast)$ (functorially). We embed $\Z$ in $\mc W$ sending $k \in \Z$ to the morphism $t \mapsto t^k$. We have an isomorphism of topological groups $\Z_p \cong \mu_{p-1} \times (1+p)^{\Z_p}$, so, if $\mc W^\circ$ is the connected component of $\mc W$ containing the trivial character, we have $\mc W = \coprod_{\mu_{p-1}^\vee} \mc W^\circ$, where $\mu_{p-1}^\vee \colonequals \Hom_{\operatorname{groups}}(\mu_{p-1}, \Q_p^\ast)$. Moreover, the map $\chi \mapsto \chi(1+p) - 1$ gives an isomorphism $\mc W^\circ \cong \mc B(1^-)$, where $\mc B(a^-)$ is the open disk of radius $a$ and center the origin. It is well known that any continuous $\chi \colon \Z_p^\ast \to A^\ast$ is locally analytic. We say that $\chi$ is \emph{$r$-admissible} if the map $\Z_p \to A^\ast$ given by $x \mapsto \chi((1+p^{r+1}x))$ is analytic. This is equivalent to
\[
|\chi(1+p) - 1| < p^{-\frac{1}{p^{r-1}(p+1)}}.
\]
Let $\chi \colon \Z_p^\ast \to K^\ast$ be $r$-admissible. It follows that there is $s \in K$ such that $\chi(t) = [t]^i \langle t \rangle^s \colonequals [t]^i \exp(s \log( \langle t \rangle)$ for all $t$ with $\val(\langle t \rangle -1) \geq r$, where
\begin{itemize}
 \item $i \in \Z/(p-1)\Z$;
 \item $[t]$ means the Teichm\"uller character applied to the reduction modulo $p$ of $t$;
 \item $\langle t \rangle t \colonequals t/[t]$ and $s \in K$ is such that $|s| < p^{-\frac{1}{p-1} + r}$.
\end{itemize}
Let $\mc W_{r^-}^\circ$ be the subset of $\mc W^\circ$ corresponding to $\mc B_{r^-} \colonequals \mc B ({p^{-\frac{1}{p^{r-1}(p+1)}}}^-)$ via $\mc W^\circ \cong \mc B(1)$. We write $\mc W_{r^-}$ for $\coprod_{\mu_{p-1}^\vee} \mc W_{r^-}^\circ$. It follows that $\chi$ is $r$-admissible if and only if $\chi \in \mc W_{r^-}(K)$. Taking \emph{closed} discs $\mc C_r$ such that $\mc B_{(r-1)^-} \subseteq \mc C_r \subseteq \mc B_{r^-}$ (where $\mc B_{0^-} \colonequals \set{0}$), we obtain an admissible open $\mc W_r \subseteq \mc W$. Furthermore $\set{\mc W_r}_{r \geq 1}$ is an admissible covering of $\mc W$ and any $\chi \in \mc W_r(K)$ is $r$-admissible (but of course the converse is not true). We fix one such $\set{\mc W_r}_{r \geq 1}$.

Let $r \geq 1$ be an integer and let $\chi \colon \Z_p^\ast \to K^\ast$ be an $r$-accessible character and write $\chi(t) = [t]^i \langle t \rangle^s \colonequals [t]^i \exp(s \log( \langle t \rangle)$, for $t$ sufficiently close to $1$. From now on, when we work at level $K(Np^r)$ we will always assume that
\[
w < \min\left( (p-1)\left (\val(s) -r - \frac{1}{p-1} \right), \frac{1}{p^{r-2}(p+1)} \right).
\]
We define the sheaf (of groups) on $\formNpn{r}(w)$
\[
\mc S_{r,w} \colonequals \Z_p^\ast (1 + p^{r-v} \mc O_{\formNpn{r}(w)}).
\]
The reason to consider $\mc S_{r,w}$ is the following: if $x \in \mc S_{r,w}(\Spf(S_r))$ (where $S_r$ is as in Section~\ref{sec: more can}) we can write $x = ub$, where $u \in \Z_p^\ast$ and $b \in 1+p^{r-v}S_r$. It follows from our assumption on $w$ that $x^\chi \colonequals \chi(u) b^s$ makes sense, and it is another section of $\mc S_{r,w}$. Thus, $\mc S_{r,w}$ is a suitable subsheaf of $\mc O_{\formNpn{r}(w)}^\ast$ where `raising to the $\chi$-th power' makes sense. To define modular forms of weight $\chi$, we need to find an analogue of the $\mc O_{\formNpn{r}(w)}^\ast$-torsor given by the generators of $\underline \omega$, but relative to $\mc S_{r,w}$, so we can raise to the $\chi$-th power. We proceed in two steps: first of all we find a suitable trivial subsheaf of $\underline \omega$ (note that, since we working with the formal models, $\underline \omega$ is not trivial even if $r \geq 1$), then, working (via the trivialization) with $\mc O_{\formNpn{r}(w)}$, we define the required torsor.

Let $\mc F_{r,w}$ be the sheaf on $\formNpn{r}(w)$ defined as follows. Let $\Spf(R) \subseteq \formN(w)$ be an open affine, and write $\Spf(S_r)$ for its pullback to $\formNpn{r}(w)$. Let $\mc G_r \colonequals \Gal(\overline R_K, S_{r,K})$. Let $\mc A$ be the false elliptic curve corresponding to $\Spf(R)$, we define
\[
\mc F_{r,w}(\Spf(S_r)) \colonequals \im(\dlog_{\mc A})^{\mc G_r} \subseteq \underline \omega_{\mc A_{S_r}/S_r}.
\]
This the desired trivial subsheaf of $\underline \omega$: by the results of Section~\ref{sec: HT} and Proposition~\ref{prop: congr E}, we have that $\mc F_{r,w} = \mc O_{\formNpn{r}(w)} E_1$. We can now define the required torsor. Writing $\mathfrak C_r^\vee$ for the (constant) sheaf defined by the dual of the canonical subgroup, we have, by Proposition~\ref{prop: HT},
\begin{gather} \label{eq: iso}
\mc F_{r,w} / p^{r-v} \mc F_{r,w} \cong \mathfrak C_r^\vee \otimes \mc O_{\formNpn{r}(w)} / p^{r-v} \mc O_{\formNpn{r}(w)}.
\end{gather}
Having this isomorphism, we define the sheaf (of sets) $\mc F'_{r,w}$ on $\formNpn{r}(w)$ as the inverse image of the constant sheaf given by ${(\mc C_r^\vee)}^\ast$, the subset of $\mc C_r^\vee$ of points of order exactly $p^r$.
\begin{prop} \label{prop: tors}
We have that $\mc F'_{r,w}$ is the trivial $\mc S_{r,w}$-torsor generated by $E_1$.
\end{prop}
\begin{proof}
This follows immediately from Section~\ref{sec: HT}.
\end{proof}
\begin{defi} \label{defi:test object p-adic 2}
A test object of level $K(Np^r)$ (with $r \geq 1$) and growth condition $w$ is a triple $((\mc A/S_r, i, \alpha, P), Y, \eta)$, where:
\begin{itemize}
 \item $S_r$ is a $p$-adically complete, flat, and normal $V$-algebra;
 \item $(\mc A / S_r, i, \alpha, P)$ is an object of the moduli problem of level$K(Np^r)$;
 \item $Y$ is a global section of $\underline \omega_{\mc A/S_r}^{\otimes 1-p}$ such that $YE_{p-1} = p^w$.
 \item $\eta$ is a global section of the pullback of $\mc F'_{r,w}$ to $\Spf(S_r)$.
\end{itemize}
\end{defi}
\begin{rmk} \label{rmk: diff test object}
The above definition is different, and not equivalent to Definition~\ref{defi:test object p-adic}, but the definition of modular form given below, in the case $\chi(t)=t^k$, will be the same as above. Since we are interested in modular forms rather than in test objects, this is not a problem.
\end{rmk}
\begin{defi} \label{defi: p-adic mod form any weight}
Let $r \geq 1$ and let $\chi$ be an $r$-accessible character. A \emph{$p$-adic modular form with respect to $D$, level $K(Np^r)$ and weight $\chi$, with coefficients in $K$ and growth condition $w$}, is a rule that assigns to every test object $T = ((\mc A/S_r, i, \alpha, P), Y, \eta)$, where $T$ is of level $K(Np^r)$ and growth condition $w$, an element $f(T) \in S_{r,K}$ such that
\begin{itemize}
 \item $f(T)$ depends only on the isomorphism class of $T$;
 \item if $\varphi \colon S_r \to S_r'$ is a morphism of $p$-adically complete, flat, and normal $V$-algebras and we denote with $T'$ the base change of $T$ to $S_r'$, we have $f(T') = \varphi(f(T))$;
 \item for any $x \in \mc S_{r,w}(\Spf(S_r)) = \Z_p^\ast(1+p^{r-v}S_r)$, we have
\[
f(((\mc A/S_r, i, \alpha, P), Y, x\eta)) = x^{\chi^{-1}} f(T).
\]
\end{itemize}
We write $S^D(K,w,K(Np^r),\chi)$ for the space of such modular forms.
\end{defi}
\begin{prop} \label{prop: p-adic mod form sheaf non int}
Let $r \geq 1$. On $\formNpn{r}(w)^{\rig}$ there is an invertible sheaf $\Omega^\chi(Np^r)(w)$ such that
\[
S^D(K,w,K(Np^r),\chi) = \Homol^0(\formNpn{r}(w)^{\rig}, \Omega^\chi(Np^r)(w)).
\]
If $\chi(t) = t^k$ for some integer $k$, we have a natural isomorphism $\Omega^\chi(Np^r)(w) \cong \Omega^k(Np^r)(w)\cong \underline \omega^{\otimes k}$, so $S^D(K,w,K(Np^r),\chi) = S^D(K,w,K(Np^r),k)$.
\end{prop}
\begin{proof}
Let $\mc O_{\formNpn{r}(w)}^{(\chi)}$ be $\mc O_{\formNpn{r}(w)}$ with the action of $\mc S_{r,w}$ twisted by $\chi$. We take for $\Omega^\chi(Np^r)(w)$ the rigidification of
\[
\shom_{\mc S_{r,w}}(\mc F'_{r,w}, \mc O_{\formNpn{r}(w)}^{(\chi^{-1})}),
\]
where $\shom_{\mc S_{r,w}}(\cdot,\cdot)$ means homomorphisms of sheaves with an action of $\mc S_{r,w}$. It is an invertible sheaf by Proposition~\ref{prop: tors}. The second part of the proposition follows by Corollary~\ref{coro: glob free}.
\end{proof}
\begin{rmk} \label{rmk: need K}
As in Remark~\ref{rmk: trivial 2}, we have a formal model of $\Omega^\chi(Np^r)(w)$. We will see that this will be false in the case $r=0$. In any case, since $E_1$ is a generator of $\underline \omega$ over $K$, we need to work after having inverted $p$ to have the isomorphism $\Omega^\chi(Np^r)(w) \cong \Omega^k(Np^r)(w)$. This is the reason why we consider modular forms with coefficients in $K$ rather than in $V$.
\end{rmk}
Let $r \geq h \geq 1$ be integers.. We write $\vartheta_{r,h} \colon \formNpn{r}(w) \to \formNpn{h}(w)$ for the natural morphism. Its rigidification is finite, \'etale, and Galois with Galois group $G_{r,h} \subseteq (\Z/p^r\Z)^\ast$ given by the image of $1+p^h\Z$ under the projection $\Z \to \Z/p^r\Z$. We let moreover $\vartheta_r \colon \formNpn{r}(w) \to \formN(w)$ be the natural morphism. We have that $\vartheta_r^{\rig}$ is finite and \'etale with $G_r = (\Z/p^r\Z)^\ast$ as Galois group. Let $\chi$ be an $r$-accessible character and let $f \in S^D(K,w,K(Np^r),\chi)$. For any $a \in G_r$, we define another modular form $f_{|\langle a \rangle}$ via the formula
\[
f_{|\langle a \rangle}(\mc A/S_r, i, \alpha, P, Y, \eta) = f(\mc A/S_r, i, \alpha, aP, Y, a^{-1}\eta),
\]
where the action of $G_r$ on the global sections of $\mc F'_{r,w}$ is defined by $aE_1 \colonequals a^\ast E_1$, viewing any $a \in G_r$ as an automorphism of $\formNpn{r}(w)$. Note that this action is compatible, via the isomorphism \eqref{eq: iso}, with the natural action of $G_r$ on ${(\mc C_r^\vee)}^\ast$. The operator $f \mapsto f_{|\langle a \rangle}$ is called a \emph{diamond operator}. It is easy to show that it comes from an action of $G_r$ on $\Omega^\chi(Np^r)(w)$.
\begin{prop} \label{prop: diamond}
Let $f \in S^D(K,w,K(Np^r),\chi)$. Then $f$ descends to level $K(Np^h)$ if and only if $f_{|\langle a \rangle} = f$ for all $a \in G_{r,h}$.
\end{prop}
\begin{proof}
This follows immediately from the definition of modular form using test objects.
\end{proof}
\begin{es}
We have $E_{1|\langle a \rangle} = [a^{-1}] E_1$ for any $a \in \m F_p^\ast$, so, as we already know, $E_1$ does not descend to level $K(N)$. On the other hand $E_{p-1|\langle a \rangle} = [a^{-1}]^{p-1} E_1^{\otimes p-1} = E_{p-1}$ as required.
\end{es}
We now study the relations between the sheaf $\Omega^\chi(Np^r)(w)$ for various $r$. Let $r \geq h \geq 1$ be integers and let $\chi$ be an $h$-accessible character. It follows that $\chi$ is also $r$-accessible.
\begin{prop} \label{prop: ind r}
We have a natural isomorphism
\[
\sigma_{r,h} \colon \left( \vartheta_{r,h,\ast} \Omega^\chi(Np^{r})(w) \right)^{G_{r,h}} \cong \Omega^\chi(Np^h)(w).
\]
Furthermore $\sigma_{r,r} = \operatorname{id}$ and, if $n \geq r$ is an integer, we have $\sigma_{n,h} = \sigma_{r,h} \circ \sigma_{n,r}$.
\end{prop}
\begin{proof}
This is \cite[Proposition~6.34]{shimura}.
\end{proof}
If $\chi$ is $r$-accessible, we have defined modular forms only of level $K(Np^r)$. We are now ready remove this restriction, and define also modular forms of level $K(Np^h)$.
\begin{defi} \label{defi: p-adic mod form any weight N}
Let $\chi \colon \Z_p \to K^\ast$ be any continuous character and let $h \geq 0$ be an integer. Choose $r\geq h$ such that $\chi$ is $r$-accessible. A \emph{$p$-adic modular form with respect to $D$, level $K(Np^h)$ and weight $\chi$, with coefficients in $K$ and growth condition $w$}, is a rule that assigns to every test object $T = ((\mc A/S_r, i, \alpha, P), Y, \eta)$, where $T$ is of level $K(Np^r)$ and growth condition $w$, an element $f(T) \in S_{r,K}$ such that the conditions of Definition~\ref{defi: p-adic mod form any weight} are satisfied and moreover
\[
f_{|\langle a \rangle} = f
\]
for all $a \in G_{r,h}$. We write $S^D(K,w,K(Np^h),\chi)$ for the space of such modular forms.
\end{defi}
\begin{rmk} \label{rmk: defi ind}
It follows from Proposition~\ref{prop: ind r}, that the above definition does not depend on the choice of $r$. Furthermore, it is obvious that Propositions~\ref{prop: diamond} and \ref{prop: ind r} stay true in general.
\end{rmk}
\begin{prop} \label{prop: mod form sheaf N}
Let $\chi$ be any continuous character and let $h$ be an integer. There is an invertible sheaf $\Omega^\chi(Np^h)(w)$ on $\formNpn{h}(w)^{\rig}$ such that
\[
S^D(K,w,K(Np^h),\chi) = \Homol^0(\formNpn{h}(w)^{\rig}, \Omega^\chi(Np^h)(w)).
\]
If $\chi(t) = t^k$ for some integer $k$, we have a natural isomorphism $\Omega^\chi(Np^h)(w) \cong \Omega^k(Np^h)(w) \cong \underline \omega^{\otimes k}$, so $S^D(K,w,K(Np^h),\chi) = S^D(K,w,K(Np^h),k)$.
\end{prop}
\begin{proof}
Choose any $r \geq h$ such that $\chi$ is $r$-accessible. We define
\[
\Omega^\chi(Np^h)(w) \colonequals \left( \vartheta_{r,h,\ast} \Omega^\chi(Np^r)(w) \right)^{G_{r,h}}.
\]
Since $\formNpn{r}(w)^{\rig} \to \formNpn{h}(w)^{\rig}$ is finite \'etale with $G_{r,h}$ as Galois group,  we have that $\Omega^\chi(Np^h)(w)$ is an invertible sheaf (see also \cite[Propositions 6.24 and 6.27]{shimura}).
\end{proof}
\begin{rmk} \label{rmk: finally}
We finally see that we do not have a formal model of $\Omega^\chi(N)(w)$. Indeed, the morphism $\formNpn{r}(w) \to \formN(w)$ is only generically \'etale.
\end{rmk}
We now study what happens for different $w$'s. Let $w' \geq w$ be a rational number that satisfies the same condition of $w$. Let $i_{w',w} \colon \formNpn{r}(w)^{\rig} \hookrightarrow \formNpn{r}(w')^{\rig}$ (for any $r$) be the natural inclusion.
\begin{prop} \label{prop: ind w}
For any continuous character $\chi$, we have a natural isomorphism
\[
\rho_{w',w} \colon i_{w',w}^\ast (\Omega^\chi(Np^r)(w')) \cong \Omega^\chi(Np^r)(w).
\]
Furthermore $\rho_{w,w} = \operatorname{id}$ and, if $w'' \geq w'$ is a rational number that satisfies the same condition of $w$, we have $\rho_{w'',w'} = \rho_{w',w} \circ \rho_{w'',w'}$.
\end{prop}
\begin{proof}
This is \cite[Lemma 6.18~]{shimura}.
\end{proof}
Let $\chi$ and $\chi'$ be continuous character and let $h$ be an integer. We say that $\chi$ and $\chi'$ are \emph{$h$-close} is the restriction of $\chi$ and $\chi'$ to $1+p^h \Z_p$ coincide. We now study the relations between modular forms of weight $\chi$ and modular forms of weight $\chi'$.
\begin{prop} \label{prop: diamond deco}
Let $0 \leq h \leq r$ be integers and let $\chi$ be a continuous character. If $\chi'$ is a continuous character that is $r$-close to $\chi$, then we have
\[
\left( \Omega^\chi(Np^r)(w)[\psi^{-1}] \right )^{G_{r,h}} = \Omega^{\chi'}(Np^h)(w),
\]
where $\psi \colonequals \chi/\chi'$ is considered as a character $G_{r,h} \to K^\ast$ and by $\Omega^\chi(Np^r)(w)[\psi^{-1}]$ we mean that the action of $G_{r,h}$ is twisted by $\psi^{-1}$.
\end{prop}
\begin{proof}
We may assume that $\chi$ and $\chi'$ are $r$-accessible. We have an isomorphism $S^D(K,w,K(Np^r),\chi) \cong S^D(K,w,K(Np^r),\chi')$ that sends a modular form $f$ to the rule $\tilde f$ defined by
\[
(\mc A/S_r, i, \alpha, P, Y, \eta) \mapsto \psi(u) f(\mc A/S_r, i, \alpha, P, Y, \eta),
\]
where $\eta = u(1+p^{r-v}x) E_1$ (recall that $\mc F'_{r,w}$ is the trivial $\mc S_{r,w}$-torsor generated by $E_1$). This isomorphism does not respect the action of the diamond operators but induces a $G_r$-equivariant isomorphism
\[
S^D(K,w,K(Np^r),\chi) \cong S^D(K,w,K(Np^r),\chi')[\psi].
\]
This means that if $f_{|\langle a \rangle a} = \psi(a)^{-1}$ for all $a \in G_{r,h}$, then $\tilde f$ is fixed under the action of $G_{r,h}$ and hence descends to level $K(Np^h)$.
\end{proof}
\begin{coro} \label{coro: deco}
Let $0 \leq h \leq r$ be integers and let $\chi$ be a continuous character. We have a $G_{r,h}$-equivariant isomorphism
\[
S^D(K,w,K(Np^r),\chi) \cong \bigoplus_{\chi ' \in I} S^D(K,w,K(Np^h),\chi'),
\]
where $I$ is a set of representative of the set of continuous characters that are $r$-close to $\chi$ modulo the equivalence relation `being $h$-close'.
\end{coro}
\begin{proof}
Let us consider the $\mc O_{\formNpn{h}(w)^{\rig}}$-module $\vartheta_{r,h,\ast}^{\rig} \Omega^\chi(Np^r)(w)$. It has an action of $G_{r,h}$ and it decompose as a direct sum of the eigenspaces of $G_{r,h}$. Hence we have a decomposition
\[
\vartheta_{r,h,\ast}^{\rig} \Omega^\chi(Np^r)(w) = \bigoplus_{\psi \colon G_{r,h} \to K^\ast} \left( \Omega^\chi(Np^r)(w)[\psi^{-1}] \right )^{G_{r,h}}.
\]
Given any $\chi' \in I$, we set $\psi \colonequals \chi/\chi'$. By definition, $\psi$ factors through a morphism $G_r \to K^\ast$ and hence it define a character $\psi \colon G_{r,h} \subseteq G_r \to K^\ast$. If $\chi'' \in I$ is different from $\chi'$, then the morphism $\varphi \colon G_{r,h} \to K^\ast$ defined by $\chi''$ must be different from $\psi$, otherwise $\chi'$ and $\chi''$ would be $h$-close. On the other hand, any $\psi \colon G_{r,h} \to K^\ast$ is obtained by some $\chi' \in I$, so the corollary follows by Proposition~\ref{prop: diamond deco}.
\end{proof}
We know explain a functorial property of the sheaves $\Omega^\chi(Np^r)(w)$ that generalizes the fact that we can pullback invariant differentials under an isogeny (see \cite[Remark~6.20]{shimura} for details). Let $r \geq 1$ be an integer and let $\chi$ be an $r$-accessible character. Let $i_{\mc A}, i_{\mc B} \colon \Spf(S_r) \to \formNpn{r}(w)$ be two morphisms, corresponding the false elliptic curves $\mc A$ and $\mc B$. Suppose we are given an isogeny $f \colon \mc B \to \mc A$ whose kernel intersects trivially the canonical subgroup of $\mc B$. We obtain an isomorphism
\[
f^\chi \colon i_{\mc B}^{\rig,\ast} \Omega^\chi(Np^r)(w) \to i_{\mc A}^{\rig,\ast} \Omega^\chi(Np^r)(w)
\]
that is the natural pullback of invariant differential if $\chi$ is an integer. If $0 \leq h \leq r$ is another integer, taking $G_{r,h}$-invariants we get a similar morphism for $\Omega^\chi(Np^h)(w)$.

Everything we have done so far can be put in analytic families over the weight space, so we see that our sheaf $\Omega^\chi(Np^r)(w)$ really `interpolate' the sheaves $\underline \omega^{\otimes k}$.
\begin{prop} \label{prop: fam ex}
For any integer $r \geq 1$ there is an invertible sheaf $\Omega(Np^h)(w)$, where $h \geq 0$, on $\mc W_r \times \formNpn{h}(w)^{\rig}$ such that if $\chi \in \mc W_r(K)$ is an $r$-accessible character, then we have a natural isomorphism
\[
(\chi, \operatorname{id})^\ast \Omega(Np^h)(w) \cong \Omega^\chi(Np^h)(w).
\]
If $r_1$ and $r_2$ are integers and $w_1$ and $w_2$ satisfy the usual condition with respect to $r_1$ and $r_2$, then the restrictions of $\Omega(Np^h)(w_1)$ and $\Omega(Np^h)(w_2)$ to $\mc W_{r_1} \cap \mc W_{r_2} \times \formNpn{h}(w_1)^{\rig} \cap \formNpn{h}(w_2)^{\rig}$ coincide. Moreover, if $0 \leq h_1 \leq h_2$, then we have a natural action of $G_{h_2,h_1}$ on $(\operatorname{id} \times \vartheta_{h_2,h_1}^{\rig})_\ast \Omega(Np^{h_2})(w)$ such that
\[
\left( (\operatorname{id} \times \vartheta_{h_2,h_1}^{\rig})_\ast \Omega(Np^{h_2})(w) \right)^{G_{h_2,h_1}} = \Omega(Np^{h_1})(w).
\]
\end{prop}
\begin{proof}
We start assuming $h \geq r$. Let $\pi_i$ be the natural projection from $\mc W_r \times \formNpn{h}(w)^{\rig}$ to the $i$-th factor. We write $\mc S_{r,v}$ for $\pi_2^{-1}\mc S_{r,w}$, and the same for $\mc F'_{r,w}$. Let $x = ub$ be a local section of $\mc S_{r,w}$. If $A \otimes B$ is a local section of $\mc O_{\mc W_r \times \formNpn{h}(w)^{\rig}}$ we define $x(A \otimes B)$ as the local section of $\mc O_{\mc W_r \times \formNpn{h}(w)^{\rig}}$ that corresponds to the function
\[
(\chi,z) \mapsto \chi(u)A(\chi)b^\chi B(z),
\]
for $\chi \in \mc W_r(T)$ and $z \in \formNpn{h}(w)^{\rig}(T)$, where $T$ is any $\Q_p$-affinoid algebra. We define $\Omega(Np^h)(w)$ as
\[
\Omega(Np^h)(w) \colonequals \shom_{\mc S_{r,w}}(\mc F'_{r,w}, \mc O_{\mc W_r \times \formNpn{h}(w)^{\rig}}).
\]
We can prove all the proposition as in the case of a single character, and we can descend to arbitrary $h$ taking $G_{r,h}$-invariants.
\end{proof}

\section{The eigencurve} \label{sec: eigen}
In this section will be mainly interested in level $K(N)$. We fix an $r$-accessible character $\chi$, where $r \geq 1$ is an integer. We have that $S^D(K,w,K(Np^r),\chi)$ is by definition the global sections of a sheaf on $\formNpn{r}(w)^{\rig}$, that is an affinoid. Furthermore, we have a formal model of $\Omega^\chi(Np^r)(w)$ that is invertible sheaf on $\formNpn{r}(w)$. It follows (see \cite[Proposition~7.1]{shimura}) that $S^D(K,w,K(Np^r),\chi)$ is a potentially orthonormizable $K$-Banach module with respect to the $\sup$ norm. For modular forms of level $K(N)$ we take $G_r$-invariants, so we have that $S^D(K,w,K(N),\chi)$ satisfies property (Pr) of \cite{buzz_eigen}, that we now recall.
\begin{defi}
Let $M$ be a Banach $A$-module, where $A$ is an affinoid $K$-algebra. Following \cite[Part I, Section 2]{buzz_eigen}, we say that $M$ satisfies the \emph{property (Pr)}, if there is a Banach $A$-module $N$ such that $M \oplus N$ is potentially orthonormizable.
\end{defi}
Let $f\in S^D(K,w,K(Np^r),\chi)$ and let $T=(\mc A/S_r,i,\alpha,P,Y,\eta)$ be a test object of level $K(Np^r)$ and growth condition $w$. There is a flat (over $V$), normal, and $p$-adically complete $S_r$-algebra $S'_r$ such that $S_{r,K} \to S_{r,K}'$ is finite \'etale and all finite and flat subgroup schemes of $\mc A[p^r]^1_K$ are defined over $S'_{r,K}$. By normality, all such subgroup schemes extend to subgroups of $\mc A[p^r]^1_{S'_r}$. We can base change $T$ to $S_r'$, obtaining the test object $T'$, and by definition of modular form we have $f(T)=f(T')$. In particular, in the following, we can assume that all finite and flat subgroup schemes of $\mc A[p^r]^1$ are defined over $S_r$. Let $\mc D$ be a finite and flat subgroup scheme of $\mc A[p^r]^1$ of order $p$ that \emph{intersect trivially the canonical subgroup of $\mc A[p^r]^1$}. By Remark~\ref{rmk: can sub}, we have that $\mc D$ defines a subgroup of $\mc A[p^r]$, and we denote with $\mc A/\mc D$ the quotient of $\mc A$ by this subgroup, that is a false elliptic curve. By our assumption on $\mc D$, the level structure on $\mc A$ defines a level structure on $\mc A/\mc D$, denoted again $(i,\alpha, P)$. Since $\mc D$ intersects trivially the canonical subgroup of $\mc A[p^r]^1$, we have an isomorphism $i_1^\ast \mc F'_{r,w} \cong i_2^\ast \mc F'_{r,w}$, where $i_1,i_2 \colon \Spf(S_r) \to \formNpn{r}(w)$ correspond to $\mc A$ and $\mc A/\mc D$, respectively. Let $\eta'$ be the image of $\eta$ under this isomorphism. By \cite[Lemma~7.5]{shimura}, we have that the $S_r$-point of $\formNpn{r}(w)$ corresponding to $(\mc A/\mc D,i,\alpha,P)$ factors through $\formNpn{r}(w/p)$. Finally, we obtain a test object $T/\mc D \colonequals ((\mc A/\mc D)/S_r,i,\alpha,Y,\eta')$ of level $K(Np^r)$ and growth condition $w/p$.
\begin{defi}
We define the $\U$-operator on $S^D(K,w,K(Np^r),\chi)$ by
\[
 \U f (T) = \frac{1}{p}\sum_{\mc D} f(T/\mc D),
\]
where the sum is over all $\mc D$ as above. Taking $G_r$-invariants, we obtain a similar operator on $S^D(K,w,K(N),\chi)$. By the above discussion, we have that the $\U$-operator factors through $S^D(K,w/p,K(N),\chi)$.
\end{defi}
Note that also the $\U$-operator can be put in families.
\begin{prop} \label{prop: U compl cont}
If $w > 0$, the $\U$ operator on $S^D(K,w,K(N),\chi)$ is completely continuous.
\end{prop}
\begin{proof}
We can work at level $K(Np^r)$. The $\U$-operator can be defined via a correspondence on $\formNpn{r}(w)$ (see \cite[Section~7]{shimura}), so it is a bounded operator. If $w >0$, the natural map $S^D(K,w/p,K(Np^r),\chi) \to S^D(K,w,K(Np^r),\chi)$ is completely continuous, so the same is true for $\U$.
\end{proof}
We have the following result of classicality for our modular forms.
\begin{prop}[{\cite[Theorem~5.1]{pay_curve}}] \label{prop: class}
Let $f \in S^D(K,w,K(N),k)$ and suppose that $\U f = a f$ for some $a\in K$. If $a$ satisfies $\val(a) < k - 1$, then $f$ is classical.
\end{prop}
Using the above notations, we now define the Hecke operators $\T_l$, for $l \neq p$ a prime that does not divide $\delta N$. Let $f\in S^D(K,w,K(Np^r),\chi)$. We may assume that all the finite and flat subgroup schemes of order $l$ of $\mc A$ are defined over $S_r$. Since in this case it is automatic that any such subgroup scheme $\mc D$ intersect trivially the canonical subgroup of $\mc A[p^r]^1$, we can define $\T_l f$ as the rule
\[
T \mapsto \frac{1}{l+1}\sum_{\mc D} f(T/\mc D).
\]
Also the operators $\T_l$ can be defined using correspondences and can be put in families. In particular, they are bounded operators, but are not completely continuous.

Let $r \geq 1$ be an integer, and assume that $w > 0$ is a rational number sufficiently small. Let $\mc Z_r$ be the spectral variety associated to the $\U$-operator acting on $\Homol^0(\mc W_r \times \formN(w)^{\rig}, \Omega(N)(w))$. We have proved that all assumptions needed to use the machine developed by Buzzard in \cite{buzz_eigen} are satisfied, so we have the following
\begin{teo} \label{teo: eigen}
We have a rigid space $\mc C_r \subseteq \mc W_r \times \m A^{1,\rig}_K$ equipped with a finite morphism $\mc C_r \to \mc Z_r$. If $L$ is a finite extension of $K$, then the points of $\mc C_r(L)$ correspond to systems of eigenvalues of modular forms with growth condition $w$, finite slope with respect to the $\U$-operator, and coefficients in $L$. If $x \in \mc C_r(L)$, let $\mc M(w)_x$ be the set of modular forms corresponding to $x$. Then all the elements of $\mc M(w)_x$ have weight $\pi_1(x) \in \mc W(L)$ and the $\U$-operator acts on $\mc M(w)_x$ with eigenvalue $\pi_2(x)^{-1}$. We have that $\pi_1$ is, locally on $\mc W_r$ and on $\mc C_r$, finite and surjective. For various $r$ and $w$, these construction are compatible. Letting $r \to \infty$ we have $w \to 0$ and we obtain the global eigencurve $\mc C \subseteq \mc W \times \m A^{1,\rig}_K$.
\end{teo}

\bibliographystyle{halpha}
\bibliography{biblio}

\begin{thebibliography}{Buz07}

\bibitem[AIS11]{over}
Fabrizio Andreatta, Adrian Iovita, and Glenn Stevens.
\newblock On overconvergent modular forms.
\newblock preprint. Available at
  \url{http://www.mathstat.concordia.ca/faculty/iovita}, 2011.

\bibitem[Bra12]{shimura}
Riccardo Brasca.
\newblock {$p$}-adic modulaf forms of non-integral weight over {S}himura
  curves.
\newblock 2012.
\newblock preprint, \href{http://arxiv.org/abs/1106.2712}{arXiv:1106.2712v4}.
  To appear in {\emph{Compositio Mathematica}}.

\bibitem[Buz07]{buzz_eigen}
Kevin Buzzard.
\newblock Eigenvarieties.
\newblock In {\em {$L$}-functions and {G}alois representations}, volume 320 of
  {\em London Math. Soc. Lecture Note Ser.}, pages 59--120. 2007.

\bibitem[Col97]{cole_class}
Robert~F. Coleman.
\newblock Classical and overconvergent modular forms of higher level.
\newblock {\em J. Th\'eor. Nombres Bordeaux}, 9(2):395--403, 1997.

\bibitem[DT94]{non-optimal}
Fred Diamond and Richard Taylor.
\newblock Nonoptimal levels of mod {$l$} modular representations.
\newblock {\em Invent. Math.}, 115(3):435--462, 1994.

\bibitem[Kas99]{pay_thesis}
Payman~L. Kassaei.
\newblock {$p$}-adic modular forms over {S}himura curves over {$\mathds Q$}.
\newblock {P}hd-{T}hesis. Available at
  \url{http://www.mth.kcl.ac.uk/~kassaei/}, 1999.

\bibitem[Kas04]{pay_totally}
Payman~L. Kassaei.
\newblock {$\mathcal P$}-adic modular forms over {S}himura curves over totally
  real fields.
\newblock {\em Compos. Math.}, 140(2):359--395, 2004.

\bibitem[Kas09]{pay_curve}
Payman~L. Kassaei.
\newblock Overconvergence and classicality: the case of curves.
\newblock {\em J. Reine Angew. Math.}, 631:109--139, 2009.

\bibitem[Lan78]{lang_cycl}
Serge Lang.
\newblock {\em Cyclotomic fields}.
\newblock Springer-Verlag, New York, 1978.
\newblock Graduate Texts in Mathematics, Vol. 59.

\bibitem[OT70]{oort}
Frans Oort and John Tate.
\newblock Group schemes of prime order.
\newblock {\em Ann. Sci. \'Ecole Norm. Sup. (4)}, 3:1--21, 1970.

\end{thebibliography}

\end{document}